\documentclass{amsart}
\usepackage{amscd,amsmath}

\newcommand{\al}{\alpha}
\newcommand{\be}{\beta}

\newcommand{\f}{\varphi}

\newcommand{\va}{\varphi}
\newcommand{\noa}{\noalign{\medskip}}

\newcommand{\hhat}[1]{\widehat{\widehat #1}}

\newcommand{\TT}{\mathbb{T}}
   
\newcommand{\RR}{\mathbb{R}}   
\newcommand{\ZZ}{\mathbb{Z}}

\numberwithin{equation}{section}

\newtheorem{te}{Theorem}[section]

\newtheorem{lm}{Lemma}[section]

\theoremstyle{definition}
\newtheorem{de}{Definition}[section]    
\newtheorem{re}{Remark}[section]
\newtheorem{ex}{Example}[section]

\begin{document}

\title{Cosphere bundle reduction in contact geometry}
\author{Oana Dr\u agulete, Liviu Ornea, Tudor S. Ratiu}
\thanks{This version: April 12, 2002}

\address{Department of Mathematics, University ``Politehnica" of
Bucharest,
\newline \indent
Bucharest, Romania}
\email{dragulete@mathem.pub.ro}

\address{University of Bucharest, Faculty of Mathematics
\newline \indent
14 Academiei str., 
70109 Bucharest, Romania}
\email{lornea@imar.ro}

\address{Institut Bernoulli, EPFL, CH-1015 Lausanne,
Switzerland}
\email{Tudor.Ratiu@epfl.ch}

\subjclass{53D20, 53D10}
\keywords{contact manifold, symplectic manifold, cotangent bundle,
cosphere bundle, momentum map, (non-zero) reduction}

\begin{abstract}
We extend the theorems concerning the equivariant symplectic reduction of the
cotangent bundle to contact geometry. The role of the cotangent bundle
is taken by the cosphere bundle. We use Albert's method for reduction
at zero and Willett's method for non-zero reduction.
\end{abstract}

\maketitle

\section{Introduction}

One of the main results concerning symplectic reduction 
with many applications in geometric mechanics states that,
in the presence of a ``good" action of a finite dimensional Lie group
$G$ on an arbitrary differentiable manifold $Q$, 
the cotangent bundle of the quotient, $T^*(Q/G)$, is symplectomorphic
with $(T^*Q)_0$, the reduced space at $0$ of the cotangent bundle. 
More generally, the reduction $T^*(Q/G)_\mu$ at $\mu\neq 0$ of $T^*Q$ 
is symplectomorphic with a vector subbundle of $T^*(Q/G_\mu)$ endowed with a 
magnetic 
symplectic form (see \cite{am}, \S 4.3; the result for $\mu = 0$ is 
due to Satzer \cite{satzer}); $G_\mu$ denotes the coadjoint isotropy subgroup at
$\mu$.

The aim of this note is to prove an analogue of this result in contact
geometry. Again we start with an arbitrary manifold $Q$ supporting a ``good"
action of a Lie group $G$. The role of the cotangent bundle will be played by 
the cosphere bundle that will be described in section 2 (cf. also \cite{ra81}). 
It
is a contact manifold. We shall prove that its reduced space at $0$ is
contactomorphic with the cosphere bundle of $Q/G$. 
Even though the result for $\mu =0$ could probably be obtained  
by ``diagram chasing", we prefer to provide an
explicit proof, identifying all contactomorphisms. More generally, we prove that 
its
reduced space at $\mu \neq 0$ embeds in a contact manner onto a subbundle of
the cosphere bundle  of $Q/G_\mu$. 

We briefly review, following \cite{al}, \cite{ge}, the reduction method at $0$ 
for 
contact manifolds. 

Recall that a \textit{contact structure} on a smooth
$(2n+1)$--dimensional manifold $N$  is a codimension one smooth distribution
$H \subset TN$, locally given by the kernel of a one-form $\eta$ such that $\eta
\wedge (d\eta)^n\ne 0$. Such an $\eta$ is called a (local) \textit{contact form}. 
Any
two proportional contact forms underly the same contact structure. 
A contact structure which is the kernel of a global contact form is
called {\em exact} or \textit{co-orientable}. If $\eta$ is 
a one form of an exact contact structure, the pair
$(N,\eta)$ is called an \textit{exact contact manifold}. 
On an exact contact manifold $N$ there is a unique vector field $R$, 
called the \textit{Reeb vector field\/}, characterized by the conditions 
$\eta (R) = 1$ and $d\eta (R, \cdot ) = 0$. The flow of the Reeb vector fields
preserves the contact form $\eta$. The Reeb vector field is nowhere vanishing
and it generates the one-dimensional distribution  $\ker d\eta =\{ v\in TN \mid
d\eta(v, \cdot)=0\}$.

 A finite dimensional
connected Lie group $G$ is said to act by \textit{contactomorphisms} on a
contact manifold if it preserves the contact structure $H$.
For an exact contact manifold $(N,\eta)$, this means that
$g^*\eta=f_g \eta$ for a smooth, real-valued, nowhere zero function $f_g$.
$G$ acts by \textit{strong contactomorphisms\/} on 
$N$, if $g^*\eta=\eta$, $i.e.$ 
$G$ preserves the contact form, not only the contact structure. 
A $G$--action by  strong
contactomorphisms on $(N, \eta)$ admits an equivariant momentum 
map $J:N\rightarrow \mathfrak{g}^*$ given by evaluating the contact 
form on fundamental fields: $\langle J, \xi\rangle = \eta(\xi_N)$. 

Throughout this paper we shall denote
by $\mathfrak{g} $ the Lie algebra of $G$, by 
$\langle \cdot, \cdot \rangle : \mathfrak{g}^\ast \times \mathfrak{g} 
\rightarrow \mathbb{R}$ the natural pairing between $\mathfrak{g}^\ast$ and  
$\mathfrak{g}$, and by $\xi_N$ the fundamental vector field (or infinitesimal 
generator) defined by $\xi \in \mathfrak{g}$. For simplicity, 
we shall work exclusively with free proper actions, although the 
extensions of our results to locally free actions is routine; in 
that case the relevant quotient spaces will be orbifolds instead 
of manifolds. For a smooth map $f: A \rightarrow B$ between the 
manifolds $A$ and $B$, $T_af :T_a A \rightarrow T_{f(a)} B$ 
denotes its derivative, or tangent map, at $a \in A$.

The momentum map $J$ is constant on the flow of the Reeb vector field. In addition, 
\[
\langle T_nJ(v), \xi \rangle = d\eta (n)(v, \xi_N(n))
\]
for any $n \in N$, $v \in T_nN$, and $\xi\in \mathfrak{g}$. This immediately implies
\[
\left[\operatorname{im}(T_nJ)\right]^\circ = \{\xi \in \mathfrak{g} \mid  
d\eta(n)(\xi_N(n), \cdot) = 0 \},
\]
which is the contact analogue of the bifurcation lemma from the usual 
theory of momentum maps on Poisson manifolds; the 
term on the left is the annihilator of the subspace in parentheses.
For this (contact) momentum map, $0\in \mathfrak{g}^*$ is  a regular value 
if and only if the fundamental fields induced by
the action do not vanish on the zero level set of $J$. Moreover,  if this is the
case, the pull back of the contact form to $J^{-1}(0)$ is basic. Let 
$\pi_0:J^{-1}(0)\rightarrow J^{-1}(0)/G$ and
$\iota_0:J^{-1}(0)\hookrightarrow N$ be the canonical projection 
and inclusion respectively. The reduction theorem asserts
the existence of a unique contact form $\eta_0$ on $J^{-1}(0)/G$
such that $\pi_0^*\eta_0=\iota_0^*\eta$. 

Regarding contact reduction at $\mu\neq 0$, up to now there are two versions 
available:
one due Albert \cite{al} and a very recent one due to Willett \cite{wil}.

\textit{Albert's method} \cite{al}. Let $(N,\eta)$ be an exact 
contact manifold with Reeb vector field
$R$ and let $\Phi$ be a ``good" action of a Lie group by strong
contactomorphisms. For $\mu \in \mathfrak{g}^\ast$, denote by 
$G_\mu$ the isotropy group at $\mu$ of the coadjoint action and by 
$\mathfrak{g}_\mu$ its Lie algebra. If $\mu\neq 0$ is a regular 
value of $J$ the restriction
of the contact form to $J^{-1}(\mu)$ is not basic. This problem is overcome by
Albert by changing the infinitesimal action of $\mathfrak{g}_\mu$ 
on $J^{-1}(\mu)$
as follows: $\xi\mapsto \xi_N- \langle\mu, \xi \rangle R$, where $R$ is the Reeb
vector field. In general, this infinitesimal action cannot be integrated to an 
action
of $G_\mu$. However, if $R$ is complete, this $\mathfrak{g}_\mu$--action is 
induced
by an action of the universal covering group $\widehat G_\mu$ (if $G_\mu$ is
connected) given by
\[
(e^{t\xi},n)\mapsto \phi_{e^{t\xi}}(\rho_{t\langle \mu,\xi \rangle}^{-1}(n)),
\]
where $\rho_t$ is the flow of the Reeb vector field. Albert defines the reduced
space as $J^{-1}(\mu)/ \widehat G_\mu$ $via$ this new action and shows it is naturally 
a contact manifold.

\textit{Willett's method} \cite{wil}. 
The idea is to expand $\mu$ and to shrink $G_\mu$. As above, $G$ is a Lie
group that acts smoothly on an exact contact manifold
$(N,\eta)$ preserving the contact form $\eta$. Let $\mu \in \mathfrak{g}^*$.  
Willett calls the \emph{kernel group of $\mu$}, the connected Lie 
subgroup $K_\mu$
of $G_\mu$  with Lie algebra $\mathfrak{k}_\mu = \ker\: 
(\mu| _{\mathfrak{g}_\mu})$. It is easy to see 
that $\mathfrak{k}_\mu$ is an ideal in $\mathfrak{g}_\mu$ 
and therefore $K_\mu$ is a connected normal subgroup of $G_\mu$.  
Contact reduction (or the contact quotient)  
of $N$ by $G$ at $\mu$ is defined by Willett as
$$
N_\mu : = J^{-1} (\RR_+ \mu)/K_\mu.
$$
Assume that $K_\mu$ 
acts freely and properly on $J^{-1} (\RR_+ \mu)$. Then $J$ is 
transversal to $\mathbb{R}_+ \mu$ and the pull back 
of $\eta$ to  $J^{-1} (\RR_+ \mu)$ is basic relative to the
$K_\mu$--action on $J^{-1} (\RR_+\mu)$ and thus induces 
a one form $\eta_\mu$ on the quotient $N_\mu$. If, in addition, 
$\ker \mu + \mathfrak{g}_\mu = \mathfrak{g}$ then the form $\eta_\mu$ 
is also a contact form. It is characterized, as usual, by the identity 
$\pi_\mu^* \eta_\mu = i_\mu^* \eta$, where  
$\pi_\mu: J^{-1} (\RR_+ \mu) \to N_\mu$ is the canonical projection 
and $i_\mu: J^{-1} (\RR_+ \mu) \hookrightarrow N$ is the canonical inclusion.

It is to be noted that for $\mu=0$, Albert's and Willett's quotients
coincide.

\medskip

\noindent {\bf Notations:} Throughout the paper we shall denote
by $\pi_G:Q\rightarrow Q/G$, $\pi_{Q/G}:T^*(Q/G)\rightarrow Q/G$,
$\pi_Q:T^*Q\rightarrow Q$ the respective canonical projections. The Liouville
one-forms of $T^*Q$ and $T^*(Q/G)$ will be denoted respectively by $\theta$ and
$\Theta$. The naturally lifted action of $G$ on $T^*Q$ admits an equivariant
momentum map $J_{ct}: T^\ast Q \rightarrow \mathfrak{g}^\ast$ given by $\langle
J_{ct}( \alpha_q), \xi \rangle = \alpha_q(\xi_Q(q))$, where $\alpha_q \in
T^\ast_qQ$, $\xi \in \mathfrak{g}$, and $\xi_Q$ denotes the fundamental vector
field defined by the $G$--action on $Q$.

\section{The cosphere bundle and its contact structure}
Let $Q$ be a differentiable manifold of real dimension $n$,
$\pi_Q:T^*Q\rightarrow Q$ its cotangent bundle, and $\theta$ the
Liouville form on $T^*Q$. We shall denote by $\al_q$, $\be_q$ etc. the
elements of $T^*Q$.

Let $G$ be a finite dimensional Lie subgroup of $\mathrm{Diff}(Q)$ and denote 
by
$\Phi:G\times Q\rightarrow Q$ a free, proper action of $G$ on $Q$. We denote by
$\Phi_*:G\times T^*Q\rightarrow T^*Q$ its natural lift to the cotangent bundle of
$Q$. $\Phi_*$ is still free and proper and preserves the Liouville form $\theta$ 
and
thus the canonical symplectic structure $-d\theta$ of $T^*Q$.

Consider  the action of the multiplicative group $\RR_+ =]0, +\infty[$ by
dilations on the fibers of  $T^*Q\setminus\{0\}$. 
\begin{de} 
The {\em  cosphere bundle} $S^\ast Q$ of $Q$ is the quotient manifold
$(T^*Q\setminus\{0\})/\RR_+$. Denote by $\kappa: [\alpha_q] \in S^\ast Q \mapsto 
q
\in Q$ the associated canonical projection.
\end{de} 
The construction described below is standard (see $e.g.$ \cite{ra81}). 

Let $\pi:T^*Q\setminus\{0\}\rightarrow S^*Q$ be the canonical
projection. The elements of the cosphere bundle are classes that we
denote with $[\al_q]$. Of course, $(\pi, \RR_+, T^*Q\setminus\{0\}, S^*Q)$ is a
$\RR_+$--principal bundle. As such, it always has global sections: it
is enough to choose a Riemannian metric on $Q$ (supposed paracompact),
to identify $T^*Q$ with $TQ$, $S^*Q$ with the unit sphere bundle
$T^1Q$ of $TQ$, and to consider the canonical inclusion $T^1Q\hookrightarrow TQ$.
Let then  $\sigma:S^*Q\rightarrow T^*Q\setminus\{0\}$ be a global section. The
equation 
\[
\sigma\circ\pi=f_\sigma 1_{T^*Q\setminus\{0\}},
\]
where $1_{T^*Q\setminus\{0\}}$ denotes the identity map of $T^*Q\setminus\{0\}$,
defines a function
$f_\sigma:T^*Q\setminus\{0\}\rightarrow \RR_+$ with the following property of
compatibility with respect to the action of
$\RR_+$:
\begin{equation}
\label{fsigma}
f_\sigma(r \al_q)=\displaystyle\frac 1r  f_\sigma(\al_q),
\quad r\in \RR_+, \; \al_q\in T^*Q\setminus\{0\}.
\end{equation} 
Indeed,
$\sigma([\al_q])=f_\sigma(\al_q)\al_q=\sigma([r\al_q])=f_\sigma(r\al_q)r\al_q.$
The following statement is now clear.
\begin{lm}
The set of global sections of $\pi:T^*Q\setminus\{0\}\rightarrow S^*Q$
is in bijective correspondence with the set of $\mathcal{C}^\infty$
functions $f:T^*Q\setminus\{0\}\rightarrow \RR_+$  satisfying \eqref{fsigma}. 
\end{lm} 
We pull back by $\sigma$ the restriction of the Liouville form and
obtain the one-form $\theta_\sigma=\sigma^*\theta$ on $S^*Q$. One has:
\begin{equation}
\label{one form relation}
\pi^*\theta_\sigma=f_\sigma\theta.
\end{equation}
Indeed, $\pi^*\theta_\sigma=\pi^*\sigma^*\theta=
(\sigma\circ\pi)^*\theta=(f_\sigma1_{T^*Q\setminus\{0\}})^*\theta
=f_\sigma\theta.$ 
Now, for another global section $\rho$, with associated function
$f_\rho$, we have
\[
\theta_\sigma=\sigma^*\theta=(\sigma\circ\pi\circ\rho)^*\theta
 =\rho^*\pi^*\theta_\sigma=\rho^*(f_\sigma\theta)=
 (f_\sigma\circ\rho)\theta_\rho,
\]
and hence we obtain
\begin{equation}\label{frho}
\theta_\sigma=g_{\sigma\rho}\theta_\rho, \quad \text{with}\quad\, 
g_{\sigma\rho}=f_\sigma\circ\rho.
\end{equation}
Note also that $g_{\sigma\rho}\circ
\pi={f_\sigma}/{f_\rho}$. From \eqref{frho} we easily
derive that \emph{$\theta_\sigma$ is a contact form on $S^*Q$ if and only if
$\theta_\rho$ is one}. But it was
proved in \cite{am} that if $\sigma$ is defined using a Riemannian
metric on $Q$, as explained above, then $\theta_\sigma$ is a contact
form. Thus we have proved:
\begin{lm} 
$\theta_\sigma$ is a global contact form on $S^*Q$ for any global
section $\sigma$. 
\end{lm}
It is also clear from \eqref{frho} that all these contact forms have
the same null space, so that the contact structure does not depend on
the choice of $\sigma$.
\begin{re}\label{con}
Let $\mathcal{C}(S^*Q)=S^*Q\times \RR_+$ be the symplectic cone over $S^*Q$,
endowed with the symplectic form $d(t\theta_\sigma)$. Then one can easily see
that $T_\sigma:\mathcal{C}(S^*Q)\rightarrow T^*Q$ given by $T_\sigma([\al_q],
t)=tf_\sigma(\al_q)\cdot\al_q$ is a well defined symplectic diffeomorphism, that 
is,
a symplectomorphism. 
\end{re}

\section{The action of $G$ on the cosphere bundle and its associated
momentum map}

We shall now lift the free proper action of $G$ to the cosphere 
bundle and compute the associated momentum  map. The action $\Phi$ 
lifts to an action  $\Phi_*$ on $T^\ast Q$ by setting 
\[
\Phi_*(g, \alpha_q):= T_{\Phi(g, q)}^\ast\Phi_{g^{-1}}\alpha_q,
\]
for $g \in G$, $\alpha_q \in T^\ast_q Q$, and where the upper star denotes
the dual map of the linear map to which it is applied. It is clear that the 
cotangent bundle projection $\pi_Q: T^\ast Q \rightarrow Q$ 
is equivariant relative to the actions $\Phi_*$ and $\Phi$. 
If the action $\Phi$ is free 
and proper, this equivariance immediately shows that the action 
$\Phi_*$ is also free and proper.

Denote by $\kappa_Q: [\alpha_q] \in S^\ast Q \mapsto q \in Q$ 
the canonical cosphere bundle projection. 

\begin{lm}\label{31}
The action $\Phi$ induces a free proper action $\widehat
\Phi_*:G\times S^*Q\rightarrow S^*Q$.
\end{lm}

\begin{proof}
Define
\[
\widehat\Phi_*(g,[\al_q])=[\Phi_*(g,\al_q)].
\]
As $\widehat\Phi_*(g,[r\al_q]) = [\Phi_*(g,r\al_q)]=[r\Phi_*(g,\al_q)]
= [\Phi_*(g,\al_q)]$, the definition is correct. Note also that $\Phi_\ast$ 
covers
$\Phi$, that is, $\kappa_Q \circ  \widehat\Phi_\ast = \Phi \circ \kappa_Q$. 
This immediately proves that freeness  (respectively 
properness) of the $G$ action  on $Q$  implies freeness 
(respectively properness) of the action 
$\widehat\Phi_\ast$ on $S^\ast Q$.
Clearly $(\widehat\Phi_{*g})^*\theta_\sigma$ is a multiple of $\theta_\sigma$ and
the proof is complete.
\end{proof}
\begin{lm}
The action $\widehat \Phi_*:G\times S^*Q\rightarrow S^*Q$ is by contactomorphisms
and the scale factors are all positive.
\end{lm}

\begin{proof}
By direct computation and using \eqref{one form relation} in the last 
equality, we have:
\begin{equation*}
\begin{split}
\widehat\Phi_{*g}^*\theta_\sigma([\alpha_q])(v_{[\al_q]})
&=\theta_{\sigma}\left(\widehat\Phi_{*g}([\al_q])\right)
\left(T_{[\alpha_q]}\widehat\Phi_{*g}(v_{[\al_ q ]})\right)\\
&= \theta\left((\sigma \circ \widehat\Phi_{*g})([\alpha_q])\right)
\left(T_{[\alpha_q]}(\sigma \circ \widehat\Phi_{*g})( v_{[\al_ q ]})\right)\\
&=  (\sigma\circ\widehat\Phi_{*g})([\al_q]) 
\left(T_{[\alpha_q]}(\pi_Q\circ\sigma\circ
\widehat\Phi_{*g})( v_{[\al_q]}) \right)\\
&= f_\sigma(\Phi_{*g}(\al_q)) \Phi_{*g}(\al_q) \left(T_{[\alpha_q]}
(\pi_Q\circ\sigma\circ \widehat\Phi_{*g})(v_{[\al_q]})\right)\\
&= f_\sigma(\Phi_{*g}(\al_q))\al_q 
\left( T_{[\alpha_q]}(\Phi_g^{-1}\circ \pi_Q\circ\sigma\circ
\widehat\Phi_{*g})(v_{[\al_q]})\right)\\
&=  f_\sigma(\Phi_{*g}(\al_q))\al_q 
\left(T_{[\alpha_q]}(\pi_Q\circ\sigma)(v_{[\al_q]})\right)\\
&= f_\sigma(\Phi_{*g}(\al_q)) \theta(\alpha_q) 
\left(T_{[\alpha_q]}(\sigma)(v_{[\al_q]})\right)\\
&= \frac{f_\sigma(\Phi_{*g}(\al_q))}{f_\sigma(\al_q)} \,\theta_\sigma([\al_q])
(v_{[\al_q]}).
\end{split}
\end{equation*}
\end{proof}

To construct a momentum  map associated to this action, we need to work with a
\emph{strong} action, that is, we need it to preserve not only the contact
structure, but the contact form. This can be achieved by adapting
Palais' argument (or, if $G$ is compact, by averaging). Indeed, owing to Lemma
\ref{31}, we may apply  Proposition 2.8 in \cite{le} asserting that for a
proper action by contactomorphisms, there always exist an invariant contact
form. (The proof of this is a straightforward modification 
of the classical proof of Palais for the existence of invariant 
Riemannian metrics on paracompact manifolds endowed with a proper 
Lie group action.) As every contact form on the cosphere bundle is 
obtained $via$ a global section as above,  we shall chose once and 
for all a section $\sigma$ for which 
$(\widehat\Phi_{*g})^*\theta_\sigma=\theta_\sigma$. Relative to this
contact form the induced action on the cosphere bundle is by strong
contactomorphisms.

The associated momentum map $J_{\theta_\sigma}$ will be denoted for simplicity
by $J$ since in what follows no other contact form different from
$\theta_\sigma$ will be used. Let
$(S^*Q)_0=J^{-1}(0)/G$ be the reduced space corresponding
to the regular value $0\in \mathfrak{g}^*$.

Similar considerations apply to the manifold $Q/G$
proving that its cosphere bundle is a contact manifold. As above, the
contact structure can be described as the kernel of a contact form of
the type $\Theta_\Sigma$, where $\Sigma:S^*(Q/G)\rightarrow
T^*(Q/G)\setminus \{0\}$ is a global section and $\Theta$ is the 
Liouville form of $T^*(Q/G)$.

\section{The main results}
We are now ready to prove:
\begin{te}\label{41}
Let $G$ be a finite dimensional Lie group, acting freely and 
properly on a differentiable manifold $Q$. Then $(S^*Q)_0$, the reduced
space at the regular value zero of the cosphere bundle of $Q$, is
contact-dif\-fe\-o\-mor\-phic  with the cosphere bundle $S^*(Q/G)$.
\end{te}

\begin{re}
Suppose $(N,\eta)$ is a contact manifold on which a Lie group $G$ acts
by strong contactomorphisms. The action can be naturally lifted to the
symplectic cone $(\mathcal{C}(N), d(t\eta))$ by letting $G$ act
trivially on $\RR_+$; one obtains an action by symplectomorphisms.
It is well known that, in this situation, the reduced symplectic space
at $0$ is the symplectic cone over the contact reduced space at $0$:
$\mathcal{C}(N_0)\cong (\mathcal{C}(N))_0$. This can be applied to
$N=S^*Q$ and combined with the cotangent bundle reduction theorem 
it should lead to a  ``diagram chasing'' proof of the theorem. However, 
we prefer to make the maps involved in the proof precise.
\end{re}

\begin{proof}
A first key observation is that the actions of $G$ and $\RR_+$ on
$T^*Q\setminus\{0\}$ commute, so that there exists the
diffeomorphism:
\begin{equation}
\label{lambda}
\lambda: (S^\ast Q)/G = (T^*Q\setminus\{0\}/\RR_+)/G \rightarrow
(T^*Q\setminus\{0\}/G)/\RR_+.
\end{equation}
Second, applying the cotangent bundle reduction theorem to $T^*Q$, we
have the symplectic diffeomorphism (see \cite{am}, \cite{ma}, or \cite{mr2})
\begin{equation}
\label{phi zero}
\varphi_0:J_{ct}^{-1}(0)/G\rightarrow T^*(Q/G),\quad \text{given~ by} \quad
 \varphi_0(\widehat{\alpha}_q)\left(T_q\pi_{G}(v_q)\right) := \al_q
(v_q),
\end{equation}
where $v_q \in T_q Q$, $\alpha_q \in J_{ct}^{-1}(0) \cap T^\ast _q Q$,
$\widehat{\alpha}_q
\in J_{ct}^{-1}(0)/G$ is its class in the reduced space at zero, and $\pi_G: Q
\rightarrow Q/G$ is the projection. Denote by $p_0: J_{ct}^{-1}(0) \rightarrow
J_{ct}^{-1}(0)/G$ the canonical projection, that is, $p_0(\alpha_q)
=\widehat{\alpha}_q$ for all $\alpha_q \in  J_{ct}^{-1}(0)$.

We want to relate the zero level sets of the contact momentum  map $J$ and 
of the symplectic momentum  map $J_{ct}$. The definition of the (contact) 
momentum map $J$, 
the $\pi$ relatedness of $\xi_{T^\ast Q}$ and 
$\xi_{S^\ast Q}$, formula \eqref{one form relation}, the definition 
of the Liouville form on $T^\ast Q$,  and finally the $\pi_Q$ 
relatedness of $\xi_{T^\ast Q}$ and $\xi_Q$ yield for any $\xi \in \mathfrak{g}$
\begin{equation*}
\begin{array}{lcl}
\langle J([\alpha_q]), \xi \rangle 
&=& \theta_\sigma([\alpha_q])\left(\xi_{S^\ast Q}([\alpha_q])\right) \\ \noa
&=& \theta_\sigma(\pi(\alpha_q))\left(T_{\alpha_q} 
\pi(\xi_{T^\ast Q}(\alpha_q)\right) \\ \noa
&=& (\pi^\ast\theta_\sigma)(\alpha_q)\left(\xi_{T^\ast Q}(\alpha_q)\right) \\ \noa
&=& f_\sigma(\alpha_q)\theta(\alpha_q)\left(\xi_{T^\ast Q}(\alpha_q)\right)\\ \noa
&=& f_\sigma(\alpha_q) \alpha_q 
\left(T_{\alpha_q} \pi_Q (\xi_{T^\ast Q}(\alpha_q)) \right) \\ \noa
&=& f_\sigma(\alpha_q) \alpha_q \left(\xi_Q (q)\right),
\end{array}
\end{equation*}
that is, 
\begin{equation}
\label{momentum map relation}
J([\alpha_q]) = f_\sigma(\alpha_q) \alpha_q.
\end{equation} 
Since $f_\sigma>0$, this implies that 
\begin{equation*}
J^{-1}(0)= \{[\al_q]\mid \al_q (\xi_Q(q))=0,
\text{~for~all~}\; \xi \in \mathfrak{g}\}.
\end{equation*}
However, $\langle J_{ct}(\alpha_q), \xi \rangle = \alpha_q(\xi_Q(q))$ 
for any $\xi \in \mathfrak{g}$, which shows that  
$J^{-1}(0)\subseteq \pi(J_{ct}^{-1}(0))$.  The
converse inclusion being obvious, we conclude that
$J^{-1}(0)=\pi(J_{ct}^{-1}(0))$ and hence
\begin{equation}
\label{zero level sets}
(S^\ast Q)_0 := J^{-1}(0)/G=(J_{ct}^{-1}(0)\setminus\{0\}/\RR_+)/G.
\end{equation}
Denote by
\begin{equation}
\label{lambda bar}
\overline{\lambda}: (S^\ast Q)_0 \rightarrow (J_{ct}^{-1}(0)
\setminus\{0\}/G)/\RR_+
\end{equation}
the diffeomorphism obtained by restricting the diffeomorphism $\lambda$ defined 
in
\eqref{lambda} to $(S^\ast Q)_0$ and denote by $\overline{\theta}_\sigma$ the
reduced contact form on $(S^\ast Q)_0$.

The definition of the diffeomorphism $\varphi_0: J_{ct}^{-1}(0)/G \rightarrow
T^*(Q/G)$ defined in \eqref{phi zero} shows that $ \widehat{0}_q \in
J_{ct}^{-1}(0)/G$ is mapped to the zero element of $T_{\pi_G(q)}^*(Q/G)$ and that
$\varphi_0$ commutes with the
$\mathbb{R}_+$--actions on $J_{ct}^{-1}(0)/G$ and on $T^*(Q/G)$ respectively. 
Thus $\varphi_0$ induces a smooth map 
\begin{equation}
\label{phi zero hat}
\widehat\varphi_0:(J_{ct}^{-1}(0)\setminus \{0\}/G)/\RR_+\rightarrow
 (T^*(Q/G)\setminus{0})/\RR_+ = S^\ast(Q/G)
\end{equation}
given by
\begin{equation}
\label{phi zero hat definition}
\widehat\f_0([\widehat{\alpha}_q]): =[\f_0(\widehat{\alpha}_q)],
\end{equation}
where $[\widehat{\alpha}_q] \in (J_{ct}^{-1}(0)/G)/\RR_+$ denotes the class of
$\widehat{\alpha}_q \in J_{ct}^{-1}(0)/G$. The same reasoning applied to
$\varphi_0^{-1}$ shows that it induces a smooth map $S^*(Q/G) \rightarrow
(J_{ct}^{-1}(0)\setminus \{0\}/G)/\RR_+$ which is easily verified to be the 
inverse
of $\widehat{\varphi}_0$, that is, $\widehat{\varphi}_0$ is a diffeomorphism.

The theorem will be proved if it is shown that $\widehat\f_0 \circ
\overline{\lambda}:(S^\ast Q)_0 \rightarrow S^\ast(Q/G)$ is a contactomorphism. 
Let $\Sigma:S^*(Q/G)\rightarrow T^*(Q/G)\setminus \{0\}$ be a global section and 
let 
$\Theta_\Sigma := \Sigma^\ast \Theta$ be the contact form on $S^*(Q/G)$ 
associated
to this section, where $\Theta$ is the Liouville form on $T^\ast(Q/G)$. From the
discussion in Section 2, we know that 
$\Theta_\Sigma$ is one of the possible contact forms underlying the contact
structure of the cosphere bundle $S^\ast(Q/G)$. Thus, to show that 
$\widehat\f_0 \circ \overline{\lambda}$ is a
contactomorphism, it will be enough to verify that
$(\widehat\f_0 \circ \overline{\lambda})^\ast \Theta_\Sigma$  is proportional to
$\overline{\theta}_\sigma$, the proportionality factor being a strictly positive
function on $(S^\ast Q)_0$. To this end, let $\pi_0:J^{-1}(0)\rightarrow
J^{-1}(0)/G = (S^\ast Q)_0$, $\iota_0:J^{-1}(0)\hookrightarrow S^*Q$ be the 
canonical
projection and the canonical inclusion, respectively. From the contact reduction
theorem at zero (reviewed in the Introduction), we know that
$\overline{\theta}_\sigma$ is characterized by the relation $\pi_0^\ast
\overline{\theta}_\sigma = \iota_0^\ast \theta_\sigma$. Thus, it suffices to show
that $(\widehat\f_0 \circ \overline{\lambda} \circ \pi_0)^\ast \Theta_\Sigma$ is
proportional to $\iota_0^\ast \theta_\sigma$ with a strictly positive function on 
$J^{-1}(0)$ as proportionality factor.

The commutative diagram below is needed in the proof that follows. All vertical
arrows are projections. The maps in this diagram have all been defined with the
exception of $\Pi: T^\ast (Q/G) \setminus \{0\} \rightarrow S^\ast(Q/G)$ which is 
the
cosphere bundle projection associated to the manifold $Q/G$ and $\overline{\pi} :
(J_{ct}^{-1}(0) \setminus \{0\})/G \rightarrow ((J_{ct}^{-1}(0) \setminus
\{0\})/G)/\mathbb{R}_+$ which is associates to each point in $(J_{ct}^{-1}(0)
\setminus \{0\})/G$ its $\mathbb{R}_+$--orbit.

\setlength{\unitlength}{1cm} 
\begin{figure}[h]
\begin{picture}(14,5.5)(2,2)
\put(2,7){$J^{-1}(0)$}
\put(4.7,7.1){\vector(-1,0){1.5}}
\put(4,7.2){$\pi$}
\put(2.5,6.7){\vector(0,-1){2.3}}
\put(2.7,5.5){$\pi_0$}
\put(2,3.7){$(S^*Q)_0$}
\put(3.3,3.8){\vector(1,0){1.5}}
\put(3.8,4){$\overline\lambda$}
\put(3.8,3.5){$\sim$}
\put(5,7){$T^*Q\supset J_{ct}^{-1}(0)\setminus\{0\}$}
\put(8.2,7){\vector(1,0){5.6}}
\put(10,7.2){$\pi_Q$}
\put(14,7){$Q$}
\put(6.9,6.7){\vector(0,-1){1}}
\put(14.1,6.7){\vector(0,-1){1}}
\put(6.4,6.2){$p_0$}
\put(13.5,6.2){$\pi_G$}
\put(5,4.5){
$
\begin{CD}
(J_{ct}^{-1}(0)\setminus\{0\})/G @>\varphi_0>\sim>T^*(Q/G)\setminus\{0\}@>\pi_{Q/G}>>Q/G\\
@V{\overline\pi}VV                                              @V{\Pi}VV\\
\left[(J_{ct}^{-1}(0)\setminus\{0\})/G\right]/\RR_+ 
@>\widehat\varphi_0>\sim> S^*(Q/G)       
\end{CD}
$
}
\end{picture}
\end{figure}

We begin with the computation of 
$(\widehat{\varphi}_0 \circ \overline{\lambda} \circ
\pi_0 \circ \pi)^\ast \Theta_\Sigma$. From the commutative diagram we have
\[
\widehat{\varphi}_0 \circ \overline{\lambda} \circ \pi_0 \circ \pi
= \Pi \circ \varphi_0 \circ p_0 \quad \text{and} \quad 
\pi_{Q/G} \circ \varphi_0 \circ p_0 = \pi_G \circ \pi_Q
\]
so that using \eqref{one form relation} with base manifold $Q/G$, 
the definition \eqref{phi zero} of
$\varphi_0$, and the global formula of the Liouville form on $T^\ast(Q/G)$, 
we get for any $\alpha_q \in J_{ct}^{-1}(0) \setminus \{0\}$ and any 
$v \in T_{\alpha_q}(J_{ct}^{-1}(0) \setminus \{0\})$
\begin{align}
\label{first pull back}
&\left((\widehat{\varphi}_0 \circ \overline{\lambda} \circ
\pi_0 \circ \pi)^\ast \Theta_\Sigma\right)(\alpha_q)(v)
= \left((\Pi \circ \varphi_0 \circ p_0)^\ast
\Theta_\Sigma\right)(\alpha_q)(v)
\nonumber \\
&\qquad = \left( (\varphi_0 \circ p_0)^\ast (\Pi^\ast \Theta_\Sigma) \right)
(\alpha_q)(v)
= \left( (\varphi_0 \circ p_0)^\ast (f_\Sigma \Theta) \right) (\alpha_q)(v)
\nonumber \\
&\qquad = (f_\Sigma \circ \varphi_0 \circ p_0)(\alpha_q) \, (\varphi_0 \circ
p_0)(\alpha_q)
\left(T_{\alpha_q}(\pi_{Q/G} \circ \varphi_0 \circ p_0)(v)\right)
\nonumber \\
&\qquad = (f_\Sigma \circ \varphi_0 \circ p_0)(\alpha_q) \, (\varphi_0 \circ
p_0)(\alpha_q)
\left(T_{\alpha_q}(\pi_G \circ \pi_Q)(v)\right)
\nonumber \\
&\qquad = (f_\Sigma \circ \varphi_0 \circ p_0)(\alpha_q) \, \alpha_q
\left(T_{\alpha_q}\pi_Q (v)\right).
\end{align}
On the other hand, since $\iota_0 :J^{-1}(0) \hookrightarrow S^\ast Q$ is the
inclusion, from  \eqref{one form relation} and the definition of the Liouville 
form
on $T^\ast Q$, we get
\begin{equation}
\label{second pull back}
\left((\iota_0 \circ \pi)^\ast \theta_\sigma \right)(\alpha_q)(v)
= (f_\sigma \theta)(\alpha_q)(v) 
= f_\sigma(\alpha_q)\, \alpha_q\left(T_{\alpha_q}\pi_Q (v)\right).
\end{equation}
The two identities \eqref{first pull back} and \eqref{second pull back} show that
on $J_{ct}^{-1}(0) \setminus \{0\}$ we have the equality 
\begin{equation}
\label{proportionality}
\pi^\ast(\widehat{\varphi}_0 \circ \overline{\lambda} \circ
\pi_0)^\ast \Theta_\Sigma =
\frac{f_\Sigma \circ \varphi_0 \circ p_0}{f_\sigma}\,
\pi^\ast \iota_0^\ast \theta_\sigma.
\end{equation}
Formula \eqref{fsigma} shows that the strictly positive proportionality factor in
\eqref{proportionality} drops to a strictly positive function $F$ on the quotient
$J^{-1}(0)$. Since $\pi$ is a surjective submersion, \eqref{proportionality} 
implies
that $(\widehat{\varphi}_0 \circ \overline{\lambda} \circ
\pi_0)^\ast \Theta_\Sigma = F \, \iota_0^\ast \theta_\sigma$ where the function
$F>0$, which is the desired identity.
\end{proof}

The first two examples below use parallelizable manifolds $Q$. Note that for an
$n$-dimensional  parallelizable manifold $Q$, the cosphere bundle is $S^\ast Q=Q\times
S^{n-1}$.

\begin{ex}\label{paral1}
Let $Q=\TT^n$ and $G=S^1$ acting by
multiplication on the first factor of the torus and trivially on the other
ones. Then $Q/G=\TT^{n-1}$ and $S^\ast(\TT^{n})=\TT^{n}\times S^{n-1}$. Hence,
by Theorem 4.1, we find that $(\TT^{n}\times S^{n-1})_0$   is contactomorphic
with   $\TT^{n-1}\times S^{n-2}$.

\end{ex}

\begin{ex}\label{paral2}
 Let $Q=\RR^n$ and $G=\ZZ^n$ acting by translations on each factor. Then
$Q/G=\TT^n$, $S^\ast(Q/G)=\TT^n\times S^{n-1}$, $S^\ast(\RR^n)=\RR^n\times
S^{n-1}$, hence we obtain the contactomorphism $ ( \RR^n\times S^{n-1})_0\cong
\TT^n\times S^{n-1}$.
\end{ex}

\begin{ex}\label{paral3}
Let $Q=S^3$ and $G=S^1$ acting by multiplication 
(of unitary quaternions by unit complex numbers). 
Then $Q/G=S^2$, the base of the Hopf fibration.  
It is well known (see, \textit{e.g.} \cite{mr1}, 
Exercise 1.2-4) that $S^\ast(S^2)$ is diffeomorphic 
with $SO(3)$. On the other hand, $S^\ast (S^3)=S^2\times S^3$. 
We thus obtain the contact diffeomorphism $( S^2\times S^3)_0\cong SO(3)$.
\end{ex}

If we want to carry out the cosphere bundle reduction at a point 
$\mu\neq 0$, we have \emph{a priori} two choices: to use
Albert's or Willett's reduction methods. 

Regarding Albert's reduction method (see its description in the
Introduction), nothing will guarantee that the action of the universal cover  
$\widehat{G}_\mu$ on $S^*Q$ is induced by an action on $Q$. Example II in
\cite{al} describes precisely such a situation. It refers  (without naming
it explicitly) to the cosphere bundle $S^*\mathbb{T}^n$ of the
$n$-dimensional torus $\mathbb{T}^n$. The group is
$G=\mathbb{T}^n$ which acts trivially on itself. For a non-zero regular value 
$\mu$ of norm $1$, Albert applies his construction with  $\widehat{G}_\mu=\RR^n$ 
and
obtains the standard circle $S^1$ as the reduced space. But the action of
$\RR^n$ on
$S^*\mathbb{T}^n$ does not come from an action of $\RR^n$ on $T^n$! Thus, 
Albert's
method cannot be used to do contact reduction of the cosphere bundle at a non 
zero value of the momentum map.

However, Willett's method can be applied, as we shall show below. 
Contact reduction at a non zero value of the momentum
map will embed in a certain cosphere bundle. The precise
statement is the following. Recall that $K_\mu$ denotes the 
connected normal Lie subgroup of $G_\mu$ whose Lie algebra is 
the ideal $\mathfrak{k}_\mu:= \ker(\mu|_{\mathfrak{g}_\mu})$ in 
$\mathfrak{g}_\mu$.

\begin{te}
Let $Q$ be a differentiable manifold of real dimension $n$, $G$ a finite 
dimensional Lie subgroup of $\operatorname{Diff}(Q)$ and 
$\Phi: G \times Q \to Q$ a smooth action  of $G$ on $Q$. 
Assume that  
$K_\mu$ acts freely and properly on $J^{-1} (\RR_+\mu)$ and that 
$\ker \mu + \mathfrak{g}_\mu = \mathfrak{g}$. Then the contact reduction  
$$
(S^*Q)_\mu = J^{-1} (\RR_+ \mu)/K_\mu
$$
is embedded by a map preserving the contact structures onto a 
subbundle of $S^* (Q/K_\mu)$.
\end{te}

\begin{proof} Consider the cosphere bundle $S^\ast Q$ endowed with the
contact form $\theta_\sigma$ preserved by the $G$--action. Willett \cite{wil}
\S 3 proves that $J$ is transversal to $\RR_+\mu$ if and only if the 
$K_\mu$--action on $J^{-1} (\RR_+\mu)$ is locally free. Our hypothesis is
that this action is in fact free, so the transversality hypothesis in 
Willett's reduction theorem is satisfied. Together with the other two
stated hypotheses, these are precisely the assumptions of Willett's 
reduction theorem reviewed in the Introduction. Thus 
$(S^*Q)_\mu = J^{-1} (\RR_+ \mu)/K_\mu$
is an exact contact manifold whose contact form, denoted by 
$\overline{\theta}_{\sigma, \mu}$, is characterized by the identity 
$\iota_\mu ^\ast \theta_\sigma = \pi_\mu^\ast \overline{\theta}_{\sigma, \mu}$,
where $\iota_\mu : J^{-1} (\RR_+\mu) \hookrightarrow S^\ast Q$
is the inclusion and 
$\pi_\mu: J^{-1} (\RR_+\mu) \rightarrow J^{-1} (\RR_+\mu)/K_\mu = (S^*Q)_\mu$ 
is the canonical projection.

As in the proof of Theorem \ref{41}, 
$J^{-1} (\RR_+ \mu)=\pi(J^{-1}_{ct}(\RR_+\mu))$ and, consequently, 
\[
(S^\ast Q)_\mu : = J^{-1}(\RR_+\mu)/K_\mu
 =(J_{ct}^{-1}(\RR_+\mu)\setminus\{0\}/\RR_+)/K_\mu.
\]
Since the actions of $K_\mu$ (by cotangent lift) 
and $\mathbb{R}_+$ (by dilation in each fiber) on $J_{ct}^{-1}(\RR_+\mu)$
commute, there is a diffeomorphism
\[
\overline{\lambda}_\mu: (S^\ast Q)_\mu \rightarrow
(J_{ct}^{-1}(\RR_+\mu)\setminus\{0\}/K_\mu)/\RR_+
\] 
characterized by the property
\[
\overline{\lambda} \circ \pi_\mu \circ \pi = \overline{\pi}_\mu \circ p_\mu,
\]
where 
\[
\overline{\pi}_\mu: (J_{ct}^{-1} (\RR_+ \mu) \setminus \{0\})/K_\mu
\rightarrow [(J_{ct}^{-1} (\RR_+ \mu) \setminus \{0\})/K_\mu]/\mathbb{R}_+
\]
and 
\[
p_\mu: J_{ct}^{-1} (\RR_+ \mu) \setminus \{0\} \subset 
T^\ast Q \setminus \{0\} \rightarrow 
(J_{ct}^{-1} (\RR_+ \mu) \setminus \{0\})/K_\mu
\]
are the canonical projections. If 
$\alpha_q \in J_{ct}^{-1} (\RR_+ \mu) \setminus \{0\}  \subset T^\ast Q$, 
denote by $\widehat\al_q = p_\mu(\alpha_q)$ its class in 
$(J_{ct}^{-1} (\RR_+ \mu) \setminus \{0\})/K_\mu$. 

Define the map
\[
\psi_\mu : (J_{ct}^{-1} (\RR_+ \mu) \setminus \{0\})/K_\mu \rightarrow 
T^* (Q/K_\mu)\setminus\{0\}
\]
by
\begin{equation}
\label{psi mu}
\psi_\mu (\widehat\al_q)(T_q \pi_{K_\mu}(v_q)) = \al_q(v_q),
\end{equation}
where $\pi_{K_\mu}: Q \rightarrow Q/K_\mu$ is the canonical projection.
To show that $\psi_\mu$ is well defined, observe that for all 
$\al_{q'} = T^*_{\phi (g,q)} \Phi_{g^{-1}} \al_q$ with $q' = \Phi (g,q)$,  
$v_{q'} = T_{\Phi (g,q)} \Phi_g (v_q + \xi_Q (q))$, and $\xi \in
\mathfrak{k}_\mu$ identity \eqref{momentum map relation} implies that
\[
\begin{array}{lcl}
\al_{q'} (v_{q'} ) &=& T^*_{\Phi(g,q)} 
\left(\Phi_{g^{-1}} \al_q, T_{\Phi(g,q)} \Phi_g (v_q + \xi_Q (q) \right) \\ \noa
&=& \al_q \left( v_q + \xi_Q (q) \right) = \al_q (v_q) + \al_q (\xi_Q(q)) \\ \noa
&=& \al_q (v_q) + 
\frac{1}{f_\sigma(\alpha_q)}\langle J([\alpha_q]), \xi \rangle \\ \noa
&=& \al_q (v_q) + 
\frac{1}{f_\sigma(\alpha_q)}\langle \mu, \xi \rangle  = \al_q (v_q) \\ \noa
\end{array} 
\]
since $\xi \in \mathfrak{k}_\mu$. This shows that $\psi_\mu$ is
well defined. It is routine to check that $\psi_\mu$ is smooth.
In addition, $\psi_\mu$ is equivariant relative to the 
$\mathbb{R}_+$--actions on $J_{ct}^{-1}(\RR_+\mu)\setminus\{0\}/K_\mu$ and
$T^\ast (Q/K_\mu) \setminus \{0\}$ respectively and 
thus it induces a smooth map on the quotients
\[
\widehat\psi_\mu:[(J_{ct}^{-1}(\RR_+\mu)
\setminus\{0\})/K_\mu]/\RR_+ \rightarrow S^*(Q/K_\mu)
\]
given by
\[
\widehat\psi_\mu([\widehat\al_q])=[\psi_\mu(\widehat \alpha_q)],
\]
where $[\widehat\al_q] := \overline{\pi}_\mu (\widehat\al_q)$, for 
$\overline{\pi}_\mu : (J_{ct}^{-1}(\RR_+\mu)
\setminus\{0\})/K_\mu \rightarrow [(J_{ct}^{-1}(\RR_+\mu)
\setminus\{0\})/K_\mu]/\RR_+$ is the canonical projection.

Next we show that $\widehat\psi_\mu$ is injective. 
If $\widehat\psi_\mu ([\widehat\al_q]) 
= \widehat\psi_\mu ([\widehat\beta_q])$, then  
there exists $r \in \RR_+$ with $\psi_\mu(\widehat\al_q) 
= r \psi_\mu(\widehat\beta_q)$, so using \eqref{psi mu}, 
$\al_q (v_q) = r \beta_q (v_q)$ for every $ v_q \in T_qQ$. 
This means that $\widehat{\alpha}_q = r \widehat{\beta}_q$
since the $K_\mu$ and $\mathbb{R}_+$ actions commute, that is,
$[\widehat{\alpha}_q] = [\widehat{\beta}_q]$ showing that  
$\widehat\psi_\mu$ is injective.

We need to show that $\widehat{\psi}_\mu \circ
\overline{\lambda}_\mu :(S^\ast Q)_\mu \rightarrow S^\ast(Q/K_\mu)$ preserves
the contact structures.  
Let $\Sigma:S^*(Q/K_\mu)\rightarrow T^*(Q/K_\mu)\setminus \{0\}$ be a 
global section and let 
$\Theta_\Sigma := \Sigma^\ast \Theta$ be the contact form on $S^*(Q/K_\mu)$ 
associated
to this section, where $\Theta$ is the Liouville form on $T^\ast(Q/K_\mu)$. The
form $\Theta_\Sigma$ is one of the possible contact forms underlying the contact
structure of the cosphere bundle $S^\ast(Q/K_\mu)$. Thus, to show that 
$\widehat\psi_\mu \circ \overline{\lambda}_\mu$ preserves the contact
structures, it will be enough to verify that
$(\widehat\psi_\mu \circ \overline{\lambda}_\mu)^\ast \Theta_\Sigma$  
is proportional to $\overline{\theta}_{\sigma, \mu}$, the 
proportionality factor being a strictly positive
function on $(S^\ast Q)_\mu$. Willett's contact reduction
theorem at $\mu \neq 0$ states that
$\overline{\theta}_{\sigma, \mu}$ is characterized by the relation 
$\pi_\mu^\ast \overline{\theta}_{\sigma, \mu} 
= \iota_\mu^\ast \theta_\sigma$. Thus, it suffices to show
that $(\widehat\psi_\mu \circ \overline{\lambda}_\mu \circ 
\pi_\mu)^\ast \Theta_\Sigma$ is
proportional to $\iota_\mu^\ast \theta_\sigma$ with a strictly 
positive function on  $J^{-1}(\mathbb{R}_+ \mu)$ as proportionality 
factor. To carry this out, we shall need a commutative diagram 
analogous to the one considered in Theorem \ref{41}.

\setlength{\unitlength}{1cm} 
\begin{figure}[h]
\begin{picture}(14,5.5)(2,2)
\put(2,7){$J^{-1}(\mathbb{R}_+ \mu)$}
\put(4.7,7.1){\vector(-1,0){1.0}}
\put(4.1,7.3){$\pi$}
\put(2.5,6.7){\vector(0,-1){2.3}}
\put(2.7,5.5){$\pi_\mu$}
\put(2,3.7){$(S^*Q)_\mu$}
\put(3.3,3.8){\vector(1,0){1.5}}
\put(3.8,4){$\overline{\lambda}_\mu$}
\put(3.8,3.5){$\sim$}
\put(5,7){$T^*Q\supset J_{ct}^{-1}(\mathbb{R}_+ \mu)\setminus\{0\}$}
\put(8.8,7.1){\vector(1,0){5.5}}
\put(10.2,7.3){$\pi_Q$}
\put(14.5,7){$Q$}
\put(6.9,6.7){\vector(0,-1){1}}
\put(14.7,6.7){\vector(0,-1){1}}
\put(6.4,6.2){$p_\mu$}
\put(14,6.2){$\pi_{K_\mu}$}
\put(5,4.5){
$
\begin{CD}
(J_{ct}^{-1}(\mathbb{R}_+ \mu)\setminus\{0\})/K_\mu @>\psi_\mu >>T^*(Q/K_\mu)\setminus\{0\}@>\pi_{Q/K_\mu}>>Q/K_\mu\\
@V{\overline\pi_\mu}VV                                              @V{\Pi_\mu}VV\\
\left[(J_{ct}^{-1}(0)\setminus\{0\})/K_\mu\right]/\RR_+ 
@>\widehat\psi_\mu>> S^*(Q/K_\mu)       
\end{CD}
$
}
\end{picture}
\end{figure}

As in the proof of Theorem \ref{41}, we begin with the computation of 
$(\widehat{\psi}_\mu \circ \overline{\lambda}_\mu \circ 
\pi_\mu \circ \pi)^\ast \Theta_\Sigma$. Since 
\[
\widehat{\psi}_\mu \circ \overline{\lambda}_\mu \circ \pi_\mu \circ \pi
= \Pi_\mu \circ \psi_\mu \circ p_\mu \quad \text{and} \quad 
\pi_{Q/K_\mu} \circ \psi_\mu \circ p_\mu = \pi_{K_\mu} \circ \pi_Q,
\]
using \eqref{one form relation} with base manifold $Q/K_\mu$, 
the definition \eqref{psi mu} of
$\psi_\mu$, and the global formula of the Liouville form on $T^\ast(Q/K_\mu)$, 
we get for any $\alpha_q \in J_{ct}^{-1}(\mathbb{R}_+ \mu) \setminus \{0\}$ 
and any $v \in T_{\alpha_q}(J_{ct}^{-1}(\mathbb{R}_+ \mu) \setminus \{0\})$
\begin{align}
\label{first mu pull back}
&\left((\widehat{\psi}_\mu \circ \overline{\lambda}_\mu \circ
\pi_\mu \circ \pi)^\ast \Theta_\Sigma\right)(\alpha_q)(v)
= \left((\Pi_\mu \circ \psi_\mu \circ p_\mu)^\ast
\Theta_\Sigma\right)(\alpha_q)(v)
\nonumber \\
&\qquad = \left( (\psi_\mu \circ p_\mu)^\ast (\Pi_\mu^\ast \Theta_\Sigma) \right)
(\alpha_q)(v)
= \left( (\psi_\mu \circ p_\mu)^\ast (f_\Sigma \Theta) \right) (\alpha_q)(v)
\nonumber \\
&\qquad = (f_\Sigma \circ \psi_\mu \circ p_\mu)(\alpha_q) \, (\psi_\mu \circ
p_\mu)(\alpha_q)
\left(T_{\alpha_q}(\pi_{Q/K_\mu} \circ \psi_\mu \circ p_\mu)(v)\right)
\nonumber \\
&\qquad = (f_\Sigma \circ \psi_\mu \circ p_\mu)(\alpha_q) \, (\psi_\mu \circ
p_\mu)(\alpha_q)
\left(T_{\alpha_q}(\pi_{K_\mu} \circ \pi_Q)(v)\right)
\nonumber \\
&\qquad = (f_\Sigma \circ \psi_0 \circ p_\mu)(\alpha_q) \, \alpha_q
\left(T_{\alpha_q}\pi_Q (v)\right).
\end{align}
On the other hand, since $\iota_\mu :J^{-1}(\mathbb{R}_+ \mu) 
\hookrightarrow S^\ast Q$ is the
inclusion, from  \eqref{one form relation} and the definition of the Liouville 
form
on $T^\ast Q$, we get
\begin{equation}
\label{second mu pull back}
\left((\iota_\mu \circ \pi)^\ast \theta_\sigma \right)(\alpha_q)(v)
= (f_\sigma \theta)(\alpha_q)(v) 
= f_\sigma(\alpha_q)\, \alpha_q\left(T_{\alpha_q}\pi_Q (v)\right).
\end{equation}
The two identities \eqref{first mu pull back} and \eqref{second mu pull back} show that
on $J_{ct}^{-1}(\mathbb{R}_+ \mu) \setminus \{0\}$ we have the equality 
\begin{equation}
\label{mu proportionality}
\pi^\ast(\widehat{\psi}_\mu \circ \overline{\lambda}_\mu \circ
\pi_\mu)^\ast \Theta_\Sigma =
\frac{f_\Sigma \circ \psi_\mu \circ p_\mu}{f_\sigma}\,
\pi^\ast \iota_\mu^\ast \theta_\sigma.
\end{equation}
Formula \eqref{fsigma} shows that the strictly 
positive proportionality factor in
\eqref{mu proportionality} drops to a strictly 
positive function $F_\mu$ on the quotient
$J^{-1}(\mathbb{R}_+ \mu)$. Since $\pi$ is a 
surjective submersion, \eqref{mu proportionality} implies
that $(\widehat{\psi}_\mu \circ \overline{\lambda}_\mu \circ
\pi_\mu)^\ast \Theta_\Sigma = F_\mu \, \iota_\mu^\ast \theta_\sigma$ 
where the function
$F_\mu >0$, which is the desired identity. This proves that 
$\widehat{\psi}_\mu \circ \overline{\lambda}_\mu$ preserves the
respective contact structures. 

That  $\widehat\psi_\mu$ is an immersion can be proved as in the embedding
version of the cotangent bundle reduction theorem (see \cite{am}, \S4.3 or
\cite{mr2}, p. 82).  Indeed,  we observe
that $\mu'=\mu|_{\mathfrak{k}_\mu}=0$, hence considering the action restricted
to $K_\mu$, the corresponding momentum map $J'$ is the restriction of $J$. We
are thus in the conditions of our Theorem 4.1 and obtain a
contact-diffeomorphism between ${J'}^{-1}(\RR_+\mu')/K_\mu 
= {J'}^{-1}(0)/K_\mu$ and $S^*(Q/K_\mu)$.
Composing this with the natural inclusion of   $J^{-1}(\RR_+\mu)/K_\mu$ in
${J'}^{-1}(\RR_+\mu')/K_\mu$, we arrive at the desired contact embedding.

This ends the proof of the theorem.
\end{proof}

 \begin{ex}
We look again at Albert's example discussed above. We have $Q=\mathbb{T}^n$, 
$G=\mathbb{T}^n$ acting naturally on itself. Take $\mu\in (\RR^n)^*$ to be 
the projection on the last factor: $\mu(x_1,\ldots,x_n)=x_n$. Then 
$\ker \mu=\RR^{n-1}$, $K_\mu=\mathbb{T}^{n-1}$ and 
$J^{-1}(\RR_+\mu)\cong \mathbb{T}^n$. Hence 
$(S^*Q)_\mu=\mathbb{T}^n/\mathbb{T}^{n-1}\cong S^1$ and $Q/K_\mu\cong S^1$. 
Our theorem yields $S^1\hookrightarrow S^*(S^1)$, the inclusion 
being the zero section in $T^*S^1$ followed by the canonical projection.
\end{ex}

\begin{ex}
We let $Q=\RR^{3(n+1)}$ and $\f$ be the natural action of $G=(\RR^3,+)$
on $Q$ by translations. The lifted action  to the cotangent bundle is
again by translations: $(\mathbf{x}, (\mathbf{q}_i,\mathbf{p}^i))\mapsto
(\mathbf{x}+\mathbf{q}_i, \mathbf{p}^i)$, $i=0,\ldots, n$. The 
symplectic momentum map (the {\em linear momentum}, see \textit{e.g.} \cite{mr2}) 
has values in $\RR^3$ (which is identified with its dual $(\RR^3)^*$ by the usual
dot product) and is given by: 
\[
J_{ct}(\mathbf{q}_i,\mathbf{p}^i)=\sum_{j=0}^n \mathbf{p}^j.
\]
Fix now $\mathbf{v}\in \RR^3\setminus \{0\}$ and define $\mu:\RR^3\rightarrow \RR$
by $ \mu(\mathbf{\xi})=\mathbf{v} \mathop{\cdot} \mathbf{\xi}$. Then we have:
\[
J_{ct}^{-1}(\RR_+\mu)=\Big\{(\mathbf{q}_i,\mathbf{p}^i) \mid
\sum_{i=0}^n \mathbf{p}^i\in \RR_+ \mathbf{v} \Big\}.
\]
As $G$ is Abelian, we have $\mathfrak{k_\mu}=\ker \mu=\mathbf{v}^\perp\cong
\RR^2$. Hence $K_\mu\cong \RR^2$. Define the map
$f: \mathbb{R}^{3(n+1)} \rightarrow \mathbb{R}^{3n+1}$ by
\[
f (\mathbf{q}_0,  \dots ,\mathbf{q}_n) := 
(\mathbf{q}_1 - \mathbf{q}_0, 
\dots, \mathbf{q}_n - \mathbf{q}_{n-1}, 
\mathbf{q}_0 \cdot \mathbf{v}/\|\mathbf{v}\|^2).
\]
Clearly $f$ is smooth, surjective, invariant under the $K_\mu$-action, and
$f(\mathbf{q}_0, \dots, \mathbf{q}_n) = f(\mathbf{q'}_0, \dots, \mathbf{q'}_n)$
if and only if  $\mathbf{q'}_i = \mathbf{q}_i + \mathbf{x}$, for all $i = 0, \dots, n$, where $\mathbf{x} \in \mathbf{v}^\perp$.  In addition, the kernel of the 
derivative of $f$ at every point equals the tangent space the $K_\mu$-orbit.
Hence $f$ induces a diffeomeorphism $Q/K_\mu = \mathbb{R}^{3(n+1)}/\mathbb{R}^2
\cong \mathbb{R}^{3n+1}$. We thus have:
\[
S^*(Q/K_\mu)=S^*
\RR^{3n+1}\cong \RR^{3n+1}\times S^{3n}.
\]
On the other hand,
\begin{equation*}
 J_{ct}^{-1} (\RR_+ \mu)/K_\mu\cong \RR^{3n+1}\times
\Big\{(\mathbf{p}^0,\ldots, \mathbf{p}^n) \mid
\sum_{i=0}^n \mathbf{p}^i \in \RR_+\mathbf{v}\Big\}
\end{equation*}
since the $K_\mu$-action does not affect $\mathbf{p}^0, \dots,
\mathbf{p}^n$.
  Applying Theorem 4.2 we have :
  $$ (S^*Q)_\mu\cong \RR^{3n+1}\times\left[ S^{3n+2}  \cap
\Big\{(\mathbf{p}^0,\ldots, \mathbf{p}^n) \mid
\sum_{i=0}^n \mathbf{p}^i \in \RR_+
\mathbf{v}\Big\}\right].$$
So Theorem 4.2 asserts the
existence of a contact structure on   the above
manifold, induced from that of $\RR^{3n+1}\times S^{3n}$. Note that it is not
obvious how to construct directly a contact structure on
$\RR^{3n+1}\times\left[ S^{3n+2}  \cap
\Big\{(\mathbf{p}^0,\ldots, \mathbf{p}^n) \mid
\sum_{i=0}^n \mathbf{p}^i \in \RR_+
\mathbf{v}\Big\}\right]$.
 \end{ex}

\begin{re}
Observe that  $\widehat\psi_\mu$ may no longer be surjective (as the 
corresponding map of the symplectic case). In fact, 
since  $\widehat\psi_\mu$ maps fibers of $J^{-1}(\RR_+\mu)/K_\mu$ in fibers of 
$S^*(Q/K_\mu)$, if it were surjective it
would be so on each fiber, but a simple count of dimensions proves this is 
impossible. On the other hand, conditions like $G=K_\mu$ or
$\mathfrak{g}=\mathfrak{g}_\mu$, which ensure surjectivity in the symplectic
case, here lead to $\mu=0$ (because of the condition $\ker \mu+
\mathfrak{g}_\mu=\mathfrak{g}$).

\end{re}

\begin{re}
There is a significant difference between the reduced spaces for $\mu \neq 0$ in 
the cotangent  bundle reduction theorem and for the cosphere bundle. The 
symplectic  quotient is  symplectically embedded (only in the particular case of 
$G$ Abelian or $G = G_\mu$ one  obtains a diffeomorphism) in $T^*(Q/G_\mu)$ 
endowed with a perturbed symplectic form (the canonical  one minus a magnetic 
term), while the contact quotient is always contactly embedded
 in  $S^*(Q/K_\mu)$
with a non perturbed contact form. Thus, in contact geometry, the cases $\mu\neq 0$ 
and $\mu = 0$ are similar and the explanation is Willett's choice of the 
kernel group of $\mu$  instead of the coadjoint isotropy group of $\mu$. Explicitly, it is 
the Lie algebra of this kernel group that 
assures the existence of a well-defined map preserving the contact structure
exactly as in the $\mu = 0$ case. 
\end{re}

\begin{re}
One may relax the assumptions on the action of $G$ (and $K_\mu$) by allowing fixed 
points and working in the category of orbifolds. 
\end{re}
\medskip

\noindent {\bf Acknowledgment.} L.O. is a member of EDGE, partially
supported by the European Contract ``Human Potential Programme Research
 Training Network" HPRN-CT-2000-00101 and by the Swiss NSF through the 
SCOPES Program during a one month
visit at the EPFL. T.S.R. was partially
supported by the European Commission and the Swiss Federal
Government through funding for the Research Training Network
\emph{Mechanics and Symmetry in Europe} (MASIE) as well as
the Swiss National Science Foundation. We would like to thank 
P. Birtea, I. Marshall, J.-P. Ortega, T. Schmah for many useful conversations 
during the elaboration of this work.

 \end{document}